\documentclass{article}

\newtheorem{fed}{\textbf{Definition}}[section]
\newtheorem{thm}[fed]{\textbf{Theorem}}
\newtheorem{lemma}[fed]{\textbf{Lemma}}

\newtheorem{prop}[fed]{\textbf{Proposition}}

\usepackage{amssymb,bbm,graphicx,epsfig,psfrag,epic,eepic,latexsym}
\usepackage{amsmath}
\usepackage{mathrsfs}
\usepackage[dvips]{color}

\begin{document}
\title{The Stark problem as a concave toric domain}
\author{Urs Frauenfelder}
\maketitle
\begin{abstract}
The Stark problem is a completely integrable system which describes the motion of an electron
in a constant electric field and subject to the attraction of a proton. In this paper we show that
in the planar case after Levi-Civita regularization the bounded component of the energy hypersurfaces of the Stark problem for energies below the critical value can be interpreted
as boundaries of concave toric domains. 
\end{abstract}

\section{Introduction}

In the Stark problem one adds to the Newtonian potential an additional linear term. This problem arises
in the study of the dynamics of an electron attracted by a proton in a constant external electric field
or a rocket attracted by a planet and subject to constant thrust.
It was
discussed already by Lagrange \cite[Section XIV]{lagrange}
as a limit case of the Euler problem of two fixed centers when one of the centers is moved to infinity.
The Stark problem is separable in parabolic coordinates as for example explained in the textbook by
Landau and Lifschitz  \cite[Section 48]{landau-lifschitz} and is therefore a completely integrable system. Its name comes from the Stark effect, namely the shifting and splitting of spectral lines in
an electric field, which was discovered independently by Stark \cite{stark} and Lo Surdo \cite{losurdo}
and whose discovery was one of the reasons that Stark received 1919 the Nobel prize. The Stark effect was one of the driving forces in the early developments of quantum mechanics. It motivated Sommerfeld to 
take account of orbits different than circular ones in quantization \cite{sommerfeld}. The separability of the Stark problem was essential in order to explain the Stark effect in the framework of Bohr-Sommerfeld quantization as was done independently by Epstein \cite{epstein} and Schwarzschild \cite{schwarzschild}.
\\ \\
In contrast to the time-honored Stark problem, the question about concave toric domains is of more recent origin. Concave toric domains play a pivotal role in symplectic embedding problems
\cite{choi-cristofaro-frenkel-hutchings-ramos, cristofaro}. Inspired by the landmark paper \cite{mcduff-schlenk} this is a very active research area in symplectic topology 
with a wealth of striking recent new results. The goal of this paper is to bring the Stark problem to the attention of this extremely prolific community by proving the following theorem.
\begin{thm}\label{main}
For every energy value below the critical one the bounded component of the regularized energy hypersurface of the planar Stark problem arises as the boundary of a concave toric domain. 
\end{thm}
We give a short overview of the paper. In Section~\ref{stark} the energy hypersurfaces of the Stark
problem are discussed. The Stark problem has a unique critical value. For energies below this critical value the energy hypersurface consists of two connected components one bounded and the other one unbounded. In this paper we restrict our attention to the bounded component for energies below the first critical value. Despite boundedness this component is not compact because of collisions. However two
body collisions can always be regularized. A classical procedure to regularize planar two body collisions
is the Levi-Civita regularisation \cite{levi-civita}. This is discussed in Section~\ref{regular}. In fact
the Levi-Civita regularization is nothing else than changing to parabolic coordinates. Therefore after applying the Levi-Civita regularization we see as well how the Stark problem separates. In particular,
the Stark problem is completely integrable and therefore foliated by Liouville tori. In Section~\ref{moment} we see that on the bounded part of the regularized energy hypersurface there are two degenerate Liouville tori having toric singularities. These correspond precisely to two collision orbits. Namely the collisions of the electron
with the proton in direction of the constant external electric force on both sides of the proton. 
Since there are only toric singularities we can define a torus action generated by a moment map. 
By a theorem of Delzant \cite[Theorem\,2.1]{delzant} the image of the moment map determines 
its preimage up to equivariant symplectomorphism. In order to prove Theorem~\ref{main} it therefore
suffices to express the image of the moment map as the graph of a convex function. This is the content
of Proposition~\ref{sc}. Since the Stark problem separates, its orbits split as the product of two
orbits belonging to perturbed harmonic oscillators. The function appearing in the proposition is
closely related to the periods of these two perturbed harmonic oscillators. Indeed, the quotient
of these periods gives the slope of the Arnold-Liouville tori.  On the other hand the periods can be
expressed with the help of the elliptic integral of the first kind. In Section~\ref{elliptic} we discuss
that the logarithm of the elliptic integral of the first kind is a strictly convex function. This
helps us in Section~\ref{strconv} to prove Proposition~\ref{sc} and hence Theorem~\ref{main}.
\\ \\
\emph{Acknowledgements: } The author acknowledges partial support by DFG grant FR 2637/2-2.

\section{The Stark problem}\label{stark}

In this section we discuss the Hamiltonian of the planar Stark problem, its energy hypersurfaces and 
their Hill's region. The potential consists of two terms, the Newtonian or Coulomb potential plus
a linear term whose derivative gives rise to a constant force, i.e.,
$$V_\varepsilon \colon \mathbb{R}^2 \setminus \{0\} \to \mathbb{R}, \quad
q \mapsto -\frac{1}{|q|}+\varepsilon q_1$$
where $\varepsilon>0$ is the field strength. The Hamiltonian is then
$$H_\varepsilon \colon T^*\big(\mathbb{R}^2 \setminus \{0\}\big) \to \mathbb{R}, \quad
(q,p) \mapsto \frac{|p|^2}{2}+V_\varepsilon(q).$$
The potential has a unique critical point
$$\mathrm{crit}(V_\varepsilon)=\big\{\big(-\tfrac{1}{\sqrt{\varepsilon}},0\big)\big\}$$
and therefore the same is true for the Hamiltonian
$$\mathrm{crit}(H_\varepsilon)=\big\{\big(-\tfrac{1}{\sqrt{\varepsilon}},0,0,0\big)\big\}.$$
The unique critical value is therefore
\begin{equation}\label{crval}
H_\varepsilon\big(-\tfrac{1}{\sqrt{\varepsilon}},0,0,0\big)=V_\varepsilon\big(-\tfrac{1}{\sqrt{\varepsilon}},0\big)=-2\sqrt{\varepsilon}.
\end{equation}
Since the Hamiltonian $H_\varepsilon$ is autonomous, i.e., independent of time, by preservation of energy it is constant along the flow lines of its Hamiltonian vector field. It seems that we have
two parameters, the energy and the field strength $\varepsilon$. However, we can get rid of one
of these parameters by a rescaling. To see that we consider for $a>0$ the diffeomorphism
$$\phi_a \colon T^*\big(\mathbb{R}^2\setminus \{0\}\big) \to T^*\big(\mathbb{R}^2\setminus \{0\}\big),
\quad (q,p) \mapsto \big(aq,\tfrac{1}{\sqrt{a}}p\big).$$
Then the Hamiltonian $H_\varepsilon$ pulls back under the diffeomorphism as
$$\phi_a^* H_\varepsilon=\frac{1}{a}H_{a^2\varepsilon},$$
the standard symplectic form
$$\omega=dp_1\wedge dq_1+dp_2\wedge dq_2$$
pulls back as
$$\phi_a^* \omega=\sqrt{a} \omega$$
so that the Hamiltonian vector field, implicitly defined by
$$dH_\varepsilon=\omega\big(\cdot,X_{H_\varepsilon}\big)$$
transforms as
$$\phi_a^* X_{H_\varepsilon}=\frac{1}{a^{3/2}}X_{H_{a^2\varepsilon}}.$$
Hence up to reparametrisation of time by a constant factor we can interpret the Hamiltonian flow
of the Stark problem for field strenght $\varepsilon$ and energy $c$ as the Hamiltonian flow 
of the Stark problem for field strength $a^2 \varepsilon$ and energy $ac$. 
\\ \\
In the following we
restrict our attention to the negative energy case. According to the previous discussion we can
assume in this case after rescaling that the energy $-\tfrac{1}{2}$. We consider therefore the
energy hypersurface
$$\Sigma_\varepsilon:=H_\varepsilon^{-1}\big(-\tfrac{1}{2}\big) \subset T^*\big(\mathbb{R}^2
\setminus \{0\}\big).$$
In view of (\ref{crval}) the energy hypersurface $\Sigma_\varepsilon$ is regular, except in the
case $\varepsilon=\tfrac{1}{16}$. If $\pi \colon T^*\mathbb{R}^2 \to \mathbb{R}^2$ denotes the footpoint
projection, we define the \emph{Hill's region} to be the shadow of the energy hypersurface 
$$\mathfrak{K}_\varepsilon:=\pi\big(\Sigma_\varepsilon\big)=\Big\{q \in \mathbb{R}^2
\setminus \{0\}: V_\varepsilon(q) \leq -\tfrac{1}{2}\Big\}.$$
The topology of the Hill's region changes dramatically at the critical field strength 
$\varepsilon=\tfrac{1}{16}$. For $\varepsilon<\tfrac{1}{16}$ the Hill's region consists of
two connected components, one bounded and the other one unbounded
$$\mathfrak{K}_\varepsilon=\mathfrak{K}_\varepsilon^b \cup \mathfrak{K}_\varepsilon^u.$$
The energy hypersurface itself decomposes into two connected components
$$\Sigma_\varepsilon=\Sigma_\varepsilon^b \cup \Sigma_\varepsilon^u$$
satisfying
$$\pi\big(\Sigma_\varepsilon^b\big)=\mathfrak{K}_\varepsilon^b, \quad
\pi \big(\Sigma_\varepsilon^u\big)=\mathfrak{K}_\varepsilon^u.$$
Note that even the component over the bounded component of the Hill's region is not itself bounded.
This is due to collision where the momenta explode. 
For $\varepsilon\geq \tfrac{1}{16}$ there just remains one unbounded component. 

\section{Regularization}\label{regular}

As we just discussed the energy hypersurface even over the bounded component is never compact due to collisions of the electron with the proton at the origin. However, two-body collisions can always be
regularized. In this section we apply a Levi-Civita regularization to the the Stark problem. An amazing 
byproduct of this regularization is, that after Levi-Civita regularization 
the Stark problem separates as well.
This shows then that the Stark problem is completely integrable. 

To describe the Levi-Civita regularization it is illuminating to use complex coordinates.
We therefore identify the configuration space $\mathbb{R}^2$ with $\mathbb{C}$ via the map
$(q_1,q_2) \mapsto q_1+iq_2$. The Levi-Civita map is defined as the two-to-one covering
$$\ell \colon \mathbb{C}\setminus \{0\} \to \mathbb{C}\setminus \{0\}, \quad z \mapsto \frac{z^2}{2}.$$
Writing down the real and imaginary part explicity, one has
$$q_1=\frac{1}{2}(z_1^2-z_2^2), \quad q_2=z_1z_2$$
which is nothing else than parabolic coordinates. 
The Levi-Civita map lifts to an exact two-to-one symplectomorphism
$$L \colon T^*\big(\mathbb{C} \setminus \{0\}\big) \to T^*\big(\mathbb{C} \setminus \{0\}\big),
\quad (z,w) \mapsto \bigg(\ell(z),\frac{w}{\overline{\ell'(z)}}\bigg)=\bigg(\frac{z^2}{2},
\frac{w}{\overline{z}}\bigg),$$
where $\overline{z}$ is the complex conjugate of $z$. 

We now define for $(z,w) \in T^*\big(\mathbb{C} \setminus \{0\}\big)$
\begin{eqnarray}\label{reg0}
E_\varepsilon(z,w)&:=&
|z|^2\Big(L^*H_\varepsilon(z,w)+\tfrac{1}{2}\Big).
\end{eqnarray}
Explicitly this becomes
\begin{equation}\label{reg}
E_\varepsilon(z,w)=E_\varepsilon^1(z_1,w_1)+E^2_\varepsilon(z_2,w_2)-2
\end{equation}
for
\begin{eqnarray*}
E_\varepsilon^1(z_1,w_1)&=&\tfrac{1}{2}w_1^2+\tfrac{1}{2}z_1^2+\tfrac{\varepsilon}{2}z_1^4\\
E_\varepsilon^2(z_2,w_2)&=&\tfrac{1}{2}w_2^2+\tfrac{1}{2}z_2^2-\tfrac{\varepsilon}{2}z_2^4.
\end{eqnarray*}
We infer two interesting consequences from formula (\ref{reg}). First we see that this formula 
makes sense for every $(z,w) \in T^*\mathbb{C}$, so that there is no
reason to take out the fiber over the origin. We therefore interpret $E_\varepsilon$ as a smooth function
$$E_\varepsilon \colon T^*\mathbb{C} \to \mathbb{R}$$
defined by (\ref{reg}). The original definition(\ref{reg0}), then leads to the equality
\begin{equation}\label{reg2}
E_\varepsilon|_{T^*(\mathbb{C}\setminus \{0\})}=R \cdot \big(L^*H_\varepsilon+\tfrac{1}{2}\big)
\end{equation}
for
$$R \colon T^* \mathbb{C} \to \mathbb{R}, \quad (z,w) \mapsto
|z|^2.$$
Adding the fiber over the origin leads to a regularization of collisions. Here are the details.
By (\ref{reg2}) we have
$$L^{-1}\big(\Sigma_\varepsilon\big)=E_{\varepsilon}|_{T^*(\mathbb{C} \setminus \{0\})}^{-1}(0)$$
and since $L$ is symplectic it further holds that
$$X_{E_\varepsilon}|_{L^{-1}(\Sigma_\varepsilon)}=R \cdot L^* X_{H_\varepsilon}|_{\Sigma_\varepsilon}.$$ 
That means that up to reparametrisation the restriction of the flow of $X_{H_\varepsilon}$
to the energy hypersurface $\Sigma_\epsilon$ can be identified with the flow of 
$X_{E_\varepsilon}$ restricted to the preimage of $\Sigma_\varepsilon$ under $L$. We now set
$$\overline{\Sigma}_\varepsilon:=E_\varepsilon^{-1}(0) \subset T^*\mathbb{C}$$
which contains $L^{-1}(\Sigma_\varepsilon)$ as a dense and open subset. The complement
$$\overline{\Sigma}_\varepsilon\setminus L^{-1}(\Sigma_\varepsilon)
=\overline{\Sigma}_\varepsilon \cap T_0\mathbb{C}$$
contains precisely the collisions where after time change the vector field now extends smoothly. 
For $0<\varepsilon<\tfrac{1}{16}$ the regularized energy hypersurface as the unregularized one decomposes into a bounded connected component and an unbounded part
$$\overline{\Sigma}_\varepsilon=\overline{\Sigma}^b_\varepsilon \cup \overline{\Sigma}^u_\varepsilon.$$
The unbounded part actually is diffeomorphic to two copies of the unbounded component
$\Sigma^u_\varepsilon$ via the two-to-one map $L$. In fact we just have
$$\overline{\Sigma}_\varepsilon^u=L^{-1}(\Sigma_\varepsilon^u).$$
The bounded part $\overline{\Sigma}_\varepsilon^b$ contains the collisions $\overline{\Sigma}_\varepsilon \cap T_0\mathbb{C}$. Different
from $\Sigma_\varepsilon^b$ it is compact. In fact it is diffeomorphic to $S^3$. This is most easily seen
by letting $\varepsilon$ go to zero, where the unbounded part disappears and the whole energy 
regularized energy hypersurface $\overline{\Sigma}_0$ becomes the standard sphere of radius $2$ in
$\mathbb{R}^4$.
\\ \\
But the Levi-Civita regularization not just regularizes the collisions but separates the problem as
well. From (\ref{reg}) we see that $E_\varepsilon$ can be written as the sum of two Poisson commuting
Hamiltonians. In particular, we see that the Stark problem is completely integrable. 
\\ \\
Due to the separability of the Stark problem we can slice our energy hypersurface. For that purpose
we abbreviate  
$$S^1_{\varepsilon,c}:=(E^1_\varepsilon)^{-1}(c) \subset T^*\mathbb{R}, 
\quad S^2_{\varepsilon,c}:=(E^2_\varepsilon)^{-1}(c) \subset T^* \mathbb{R}.$$
Note that 
$$\mathrm{crit}E^2_\varepsilon=\big\{\big(0,0\big), \big(\tfrac{1}{\sqrt{2\varepsilon}},0\big),
\big(-\tfrac{1}{\sqrt{2\varepsilon}},0\big)\big\}.$$
The first critical point is a local minimum and the two other critical points are saddle points.
Its critical values are 
$$E^2_\varepsilon(0,0)=0, \quad  E^2_\varepsilon\big(\pm \tfrac{1}{\sqrt{2\varepsilon}},0\big)=\tfrac{1}{8\varepsilon}.$$
For $0\leq c <\tfrac{1}{8\varepsilon}$ the set $S^2_{\varepsilon,c}$ decomposes as
$$S^2_{\varepsilon,c}=S^{2,b}_{\varepsilon,c} \cup S^{2,u}_{\varepsilon,c}$$
where $S^{2,b}_{\varepsilon,c}$ is the bounded component and the unbounded part
$S^{2,u}_{\varepsilon,c}$ consists of two connected components which are symmetric to each other
via the reflection $(z_2,w_2) \mapsto (-z_2,w_2)$. If $0<c<\tfrac{1}{8\varepsilon}$ the bounded component
$S^{2,b}_{\varepsilon,c}$ is diffeomorphic to a circle. For $c=0$ it degenerates to a point. Different
from $S^2_{\varepsilon,c}$ the set $S^1_{\varepsilon,c}$ is always bounded. For positive $c$ it
is diffeomorphic to a circle, for $c=0$ it degenerates to a point and for negative $c$ it is empty.
For $0<\varepsilon<\tfrac{1}{16}$ we have the slicing 
\begin{equation}\label{slice}
\overline{\Sigma}_c^b=\bigcup_{0 \leq c \leq 2} S^1_{\varepsilon,2-c} \times
S^{2,b}_{\varepsilon,c}.
\end{equation}
If $0<c<2$, then $S^1_{\varepsilon,2-c} \times
S^{2,b}_{\varepsilon,c}$ is a torus, namely an Arnold-Liouville torus expected for a completely integrable system. For $c=0$ or $c=2$ the Arnold-Liouville torus degenerates to a circle. These circles
are then periodic orbits. Going back to the unregularized system they correspond to collision orbits. 
For $c=0$ it is the collision orbit on the positive $q_1$-ray and for $c=2$ it is the one on 
the negative $q_1$-ray. 

\section{The moment map}\label{moment}

In this section we assume that the field strength satisfies
$$0<\varepsilon<\frac{1}{16}.$$
We first define a torus action on the regularized moduli space $\overline{\Sigma}_\varepsilon^b$. In order to do that we first need the periods. For $c>0$ the set $S^1_{\varepsilon,c}$ is diffeomorphic
to a circle, which coincides with the periodic orbit of the Hamiltonian $E^1_\varepsilon$ of
energy $c$. By Hamiltons equation we have $\dot{z}_1=w_1$ so that by definition of $E^1_\varepsilon$ we obtain
\begin{equation}\label{e1}
\frac{\dot{z}^2_1}{2}+\frac{z_1^2}{2}+\frac{\varepsilon z_1^4}{2}=c.
\end{equation}
We parametrize this periodic orbit such that it starts at time zero with zero velocity at its maximum.
From (\ref{e1}) we see that at its maximum $z_1^2$ satisfies the quadratic equation
$$\varepsilon z_1^4+z_1^2-2c=0$$
so that
$$z_1^2=\frac{-1+\sqrt{1+8c \varepsilon}}{2\varepsilon}.$$
It then gets accelerated to the centre such that after some time it passes the origin. Then it decelerates symmetrically with respect to the origin such that after the same amount of time
it attains the minimum. After that we can let the movie run backwards until it attains the maximum again.
The time it takes from the maximum to the origin is therefore precisely a quarter of the period. 
With the help of (\ref{e1}) we compute this quarter period as
\begin{eqnarray*}
\frac{\tau^1_\varepsilon(c)}{4}&=&\int_0^{\frac{\tau^1_\varepsilon(c)}{4}}dt\\
&=&\int_0^{\sqrt{\frac{-1+\sqrt{1+8c \varepsilon}}{2\varepsilon}}}\frac{1}{\sqrt{2c-z_1^2-\varepsilon z_1^4}}dz_1\\
&=&\frac{1}{\sqrt{\varepsilon}}\int_0^{\sqrt{\frac{-1+\sqrt{1+8c \varepsilon}}{2\varepsilon}}}\frac{1}{\sqrt{\Big(\frac{-1+\sqrt{1+8c \varepsilon}}{2\varepsilon}-z_1^2\Big)\Big(\frac{1+\sqrt{1+8c \varepsilon}}{2\varepsilon}+z_1^2\Big)}}dz_1.
\end{eqnarray*}
Using the change of variables 
$$\zeta=\frac{z_1}{\sqrt{\frac{-1+\sqrt{1+8c \varepsilon}}{2\varepsilon}}}$$
we rephrase this as 
\begin{eqnarray}\label{per1}
\tau^1_\varepsilon(c)&=&\frac{4}{\sqrt{\varepsilon}}\int_0^1 \frac{1}{\sqrt{\Big(1-\zeta^2\Big)
\Big(\frac{1+\sqrt{1+8c\varepsilon}}{2\varepsilon}+\frac{-1+\sqrt{1+8c\varepsilon}}{2\varepsilon}\zeta^2\Big)}}d\zeta\\ \nonumber
&=&\frac{2^{5/2}}{\sqrt{1+\sqrt{1+8c\varepsilon}}}\int_0^1 \frac{1}{\sqrt{\Big(1-\zeta^2\Big)
\Big(1-\frac{1-\sqrt{1+8c\varepsilon}}{1+\sqrt{1+8c\varepsilon}}\zeta^2\Big)}}d\zeta\\ \nonumber
&=&\frac{2^{5/2}}{\sqrt{1+\sqrt{1+8c\varepsilon}}}K\bigg(\frac{1-\sqrt{1+8c\varepsilon}}{1+\sqrt{1+8c\varepsilon}}\bigg)
\end{eqnarray}
where
$$K(m):=\int_0^1 \frac{1}{\sqrt{1-\zeta^2)(1-m\zeta^2)}}d\zeta,
\quad m \in (-\infty,1)$$
is the elliptic integral of the first kind.
\\ \\
Similarly for $0<c<\tfrac{1}{8\varepsilon}$ the set $S^{2,b}_{\varepsilon,c}$ corresponds to the
trace of a periodic orbit satisfying
$$\frac{\dot{z}^2_2}{2}+\frac{z_2^2}{2}-\frac{\varepsilon z_2^4}{2}=c.$$
Replacing in the above computation $\varepsilon$ by $-\varepsilon$ we obtain for its period
\begin{equation}\label{per2}
\tau^2_\varepsilon(c)=\frac{2^{5/2}}{\sqrt{1+\sqrt{1-8c\varepsilon}}}K\bigg(\frac{1-\sqrt{1-8c\varepsilon}}{1+\sqrt{1-8c\varepsilon}}\bigg).
\end{equation}
Denote by $\phi^t_{E_\varepsilon^1}$ the flow of the Hamiltonian vector field of $E_\varepsilon^1$
on $T^*\mathbb{R}$ and by $\phi^t_{E_\varepsilon^2}$ the flow of the Hamiltonian vector field
of $E_\varepsilon^2$. We abbreviate by $S^1=\mathbb{R}/\mathbb{Z}$ the circle and define the
two-dimensional torus as $T^2=S^1 \times S^1$. In view of the slicing (\ref{slice}) we are now in position to define a torus action
$$T^2 \times \overline{\Sigma}_\varepsilon^b \to \overline{\Sigma}_\varepsilon^b$$
given by
$$(t_1,t_2,z_1,w_1,z_2,w_2) \mapsto 
\bigg(\phi^{t_1\tau^1_\varepsilon(E^1_\varepsilon(z_1,w_1))}_{E^1_\varepsilon}(z_1,w_1),
\phi^{t_2\tau^2_\varepsilon (E^2_\varepsilon(z_2,w_2))}_{E^2_\varepsilon}(z_2,w_2)\bigg).$$
Let $\mathcal{T}^1_\varepsilon$ be the primitive of $\tau^1_\varepsilon$ given by
$$\mathcal{T}^1_\varepsilon(c)=\int_0^c \tau^1_\varepsilon(b) db$$
and similarly define
$$\mathcal{T}^2_\varepsilon(c)=\int_0^c \tau^2_\varepsilon(b) db.$$
Then the map
$$\mu_\varepsilon=(\mu^1_\varepsilon,\mu^2_\varepsilon) \colon \overline{\Sigma}_\varepsilon^b \to \mathbb{R}^2=\mathrm{Lie}(T^2)$$
with
$$\mu^1_\varepsilon=\mathcal{T}^1_\varepsilon \circ E_\varepsilon^1, \quad \mu^2_\varepsilon=\mathcal{T}^2_\varepsilon \circ E_\varepsilon^2$$
is a moment map for the torus action on $\overline{\Sigma}_\varepsilon^b$. By the slicing (\ref{slice})
its image is given by
$$\mathrm{im}\mu_\varepsilon=\Big\{\big(\mathcal{T}^1_\varepsilon(2-c),\mathcal{T}^2_\varepsilon(c)\big): c \in [0,2]\Big\}
\subset [0,\infty)^2 \subset \mathbb{R}^2.$$
The functions $\mathcal{T}_\varepsilon^1$ and $\mathcal{T}_\varepsilon^2$ are both strictly monotone.
Therefore there exists a strictly decreasing smooth function
$$f_\varepsilon \colon \big[0,\mathcal{T}_\varepsilon^1(2)\big] \to \big[0,\mathcal{T}_\varepsilon^2(2)\big]$$
such that 
\begin{equation}\label{gr}
\mathcal{T}^2_\varepsilon(c)=f_\varepsilon\big(\mathcal{T}_\varepsilon^1(2-c)\big), 
\quad c \in [0,2].
\end{equation}
Note that the image of the moment map can be written as the graph
$$\mathrm{im}\mu_\varepsilon=\Gamma_{f_\varepsilon}.$$
Since by Delzant \cite{delzant} the image of the moment map determines its preimage up to
equivariant symplectomorphisms Theorem~\ref{main} follows from the following proposition. 
\begin{prop}\label{sc}
For any $0<\varepsilon<\tfrac{1}{16}$ it holds that
$$f_\varepsilon''(x)>0, \quad x \in \big[0,\mathcal{T}_\varepsilon^1(2)\big],$$
i.e., the function $f_\varepsilon$ is strictly convex. 
\end{prop}
We prove the Proposition in Section~\ref{strconv}. In this section we just want to express
the second derivative of $f_\varepsilon$ in terms of the period functions and its logarithmic
derivatives. 
\\ \\
Differentiating (\ref{gr}) we obtain for any $c \in [0,2]$
\begin{eqnarray*}
\tau^2_\varepsilon(c)&=&\big(\mathcal{T}^2_\varepsilon\big)'(c)\\
&=&-f_\varepsilon'\big(\mathcal{T}^1_\varepsilon(2-c)\big)\big(\mathcal{T}^1_\varepsilon\big)'(2-c)\\
&=&-f_\varepsilon'\big(\mathcal{T}^1_\varepsilon(2-c)\big)\tau^1_\varepsilon(2-c),
\end{eqnarray*}
which we rewrite as
$$f_\varepsilon'\big(\mathcal{T}^1_\varepsilon(2-c)\big)=-\frac{\tau^2_\varepsilon(c)}{\tau^1_\varepsilon(2-c)}.$$
Differentiating this equality once more we obtain
\begin{eqnarray*}
-f_\varepsilon''\big(\mathcal{T}^1_\varepsilon(2-c)\big)\tau^1_\varepsilon(2-c)
&=&-\frac{(\tau^2_\varepsilon)'(c)\cdot \tau_\varepsilon^1(2-c)+\tau_\varepsilon^2(c)
\cdot(\tau_\varepsilon^1)'(2-c)}{(\tau_\varepsilon^1(2-c))^2}
\end{eqnarray*}
and therefore
\begin{equation}\label{logper}
f''_\varepsilon\big(\mathcal{T}_\varepsilon^1(2-c)\big)=\frac{\tau^2_\varepsilon(c)}
{(\tau_\varepsilon^1(2-c))^2}\Big((\ln \tau^2_\varepsilon)'(c)+(\ln \tau^1_\varepsilon)'(2-c)\Big).
\end{equation}
In Section~\ref{strconv} we use this formula to prove Proposition~\ref{strconv}.

\section{Elliptic integrals}\label{elliptic}

Recall that the elliptic integral of the first kind is defined as 
$$K(m):=\int_0^1 \frac{1}{\sqrt{(1-\zeta^2)(1-m \zeta^2)}}d\zeta,\quad  m \in (-\infty,1).$$
In order to prove strict convexity in Proposition~\ref{sc} we need the following lemma.
\begin{lemma}\label{moninc}
For every $m \in (-\infty,1)$ it holds that
$$(\ln K)''(m)>0,$$
i.e. the logarithm of $K$ is strictly convex. 
\end{lemma}
\textbf{Proof: } The first two derivates of $K$ are given by
$$K'(m)=\frac{1}{2}\int_0^1 \frac{\zeta^2}{\sqrt{(1-\zeta^2)(1-m \zeta^2)^3}}d\zeta$$
$$K''(m)=\frac{3}{4}\int_0^1 \frac{\zeta^4}{\sqrt{(1-\zeta^2)(1-m \zeta^2)^5}}d\zeta$$
Using Cauchy-Schwarz inequality we obtain the following interpolation inequality
\begin{eqnarray*}
K'(m)&=&\frac{1}{2}\int_0^1 \frac{\zeta^2}{\sqrt{(1-\zeta^2)(1-m \zeta^2)^3}}d\zeta\\
&=&\frac{1}{2}\int_0^1 \frac{1}{(1-\zeta^2)^{1/4}(1-m\zeta^2)^{1/4}}
\cdot \frac{\zeta^2}{(1-\zeta^2)^{1/4}(1-m\zeta^2)^{5/4}}d\zeta\\
&\leq&\frac{1}{2}\sqrt{\int_0^1 \frac{1}{\sqrt{(1-\zeta^2)(1-m \zeta^2)}}d\zeta}
\cdot \sqrt{\int_0^1 \frac{\zeta^4}{\sqrt{(1-\zeta^2)(1-m \zeta^2)^5}}d\zeta}\\
&=&\frac{1}{2}\sqrt{K(m)}\cdot \sqrt{\frac{4}{3}K''(m)}\\
&=&\frac{1}{\sqrt{3}}\sqrt{K(m) \cdot K''(m)}
\end{eqnarray*}
and therefore
$$K(m) K''(m) \geq 3\big(K'(m)\big)^2.$$
We infer from this for the derivative of the logarithmic derivative of $K$
\begin{eqnarray*}
\bigg(\frac{K'(m)}{K(m)}\bigg)'&=&\frac{K''(m)K(m)-\big(K'(m)\big)^2}{\big(K(m)\big)^2}\\
&\geq&\frac{3\big(K'(m)\big)^2-\big(K'(m)\big)^2}{\big(K(m)\big)^2}\\
&=&2\bigg(\frac{K'(m)}{K(m)}\bigg)^2\\
&>&0
\end{eqnarray*}
This proves the lemma. \hfill $\square$

\section{Proof of strict convexity}\label{strconv}

In this section we prove Proposition~\ref{strconv} and therefore Theorem~\ref{main}. For that purpose
we introduce the function
$$\Phi \colon (-\infty,1) \to \mathbb{R}, \quad
x \mapsto \frac{1}{\sqrt{1+\sqrt{1-x}}}K\bigg(\frac{1-\sqrt{1-x}}{1+\sqrt{1-x}}\bigg).$$
Its derivative is given by
\begin{eqnarray*}
\Phi'(x)&=&\frac{1}{4\sqrt{(1-x)(1+\sqrt{1-x})^3}}K\bigg(\frac{1-\sqrt{1-x}}{1+\sqrt{1-x}}\bigg)\\
& &+\frac{1}{\sqrt{(1-x)(1+\sqrt{1-x})^5}}K'\bigg(\frac{1-\sqrt{1-x}}{1+\sqrt{1-x}}\bigg)
\end{eqnarray*}
and therefore its logarithmic derivative reads
\begin{eqnarray*}
(\ln \Phi)'(x)&=&\frac{1}{4\sqrt{1-x}(1+\sqrt{1-x})}\\
& &+\frac{1}{\sqrt{1-x}(1+\sqrt{1-x})^2}
(\ln K)'\bigg(\frac{1-\sqrt{1-x}}{1+\sqrt{1-x}}\bigg).
\end{eqnarray*}
The functions $x \mapsto \tfrac{1-\sqrt{1-x}}{1+\sqrt{1-x}}$, $x \mapsto \frac{1}{4\sqrt{1-x}(1+\sqrt{1-x})}$, and $x \mapsto \frac{1}{\sqrt{1-x}(1+\sqrt{1-x})^2}$ are all strictly monotone increasing, therefore by Lemma~\ref{moninc} we conclude that the same is true for the logarithmic derivative
of $\Phi$. 
\\ \\
From (\ref{per1}) and (\ref{per2}) we see that
$$\tau^1_\varepsilon(c)=2^{5/2} \Phi(-8\varepsilon c), \quad
\tau^2_\varepsilon(c)=2^{5/2} \Phi(8\varepsilon c)$$
and hence
\begin{eqnarray*}
(\ln \tau_\varepsilon^2)'(c)+(\ln \tau_\varepsilon^1)'(2-c)&=&8\varepsilon
\Big((\ln \Phi)'(8\varepsilon c)-(\ln \Phi)'(8\varepsilon c-16 \varepsilon)\Big)>0
\end{eqnarray*}
where we used for the inequality that the logarithmic derivative of $\Phi$ is strictly monotone increasing. Combined with (\ref{logper}) this implies that
$$f_\varepsilon''(x)>0$$
for every $x \in [0,\mathcal{T}_\varepsilon^1(2)]$. This finishes the proof of Proposition~\ref{strconv}
and of Theorem~\ref{main}.

\end{document}